\documentclass[11pt]{amsart}

\usepackage{amsmath, amsthm, amssymb, amsfonts, mathrsfs, amscd, amsthm}
\usepackage[all]{xy}

\newtheorem{thm}{Theorem}[section]
\newtheorem{pro}[thm]{Proposition}
\newtheorem{lem}[thm]{Lemma}
\newtheorem{cor}[thm]{Corollary}

\theoremstyle{definition}
\newtheorem{de}[thm]{Definition}

\theoremstyle{remark}

\newcommand{\mc}[1]{\mathcal{#1}}
\newcommand{\mb}[1]{\mathbb{#1}}
\newcommand{\mf}[1]{\mathfrak{#1}}

\newcommand{\defi}{\begin{de}}
\newcommand{\ed}{\end{de}}
\newcommand{\bl}{\begin{lem}}
\newcommand{\el}{\end{lem}}
\newcommand{\bp}{\begin{pro}}
\newcommand{\ep}{\end{pro}}
\newcommand{\bt}{\begin{thm}}
\newcommand{\et}{\end{thm}}
\newcommand{\bc}{\begin{cor}}
\newcommand{\ec}{\end{cor}}
\newcommand{\bpf}{\begin{proof}}
\newcommand{\epf}{\end{proof}}

\newcommand{\beq}{\begin{equation}}
\newcommand{\eeq}{\end{equation}}
\newcommand{\beqs}{\begin{equation*}}
\newcommand{\eeqs}{\end{equation*}}
\newcommand{\ben}{\begin{enumerate}}
\newcommand{\een}{\end{enumerate}}
\newcommand{\bit}{\begin{itemize}}
\newcommand{\eit}{\end{itemize}}

\begin {document}

\title{Very ample polarized self maps extend to projective space}

\author{Anupam Bhatnagar}
\address{Anupam Bhatnagar; Department of Mathematics and Computer Science; Lehman College;
250 Bedford Park Boulevard West, Bronx, NY 10468 U.S.A.}
\email{anupambhatnagar@gmail.com}

\author{Lucien Szpiro}
\address{Lucien Szpiro; Ph.D. Program in Mathematics; CUNY Graduate Center; 365 Fifth Avenue, New York, NY 10016-4309 U.S.A.}
\email{lszpiro@gc.cuny.edu}

\date{June 15, 2011. \emph{Keywords and Phrases:} Arithmetic Dynamical Systems on Algebraic Varieties.
 Both authors are partially supported by NSF Grants DMS-0854746 and DMS-0739346.}

\begin{abstract}
Let $X$ be a projective variety defined over an infinite field, equipped with a line bundle $L$, giving an embedding of $X$ into $\mb{P}^m$ and let $\phi: X \to X$ be a morphism such that $\phi^*L \cong L^{\otimes q}, q\geq 2$. Then there exists an integer $r>0$ extending $\phi^r$ to $\mb{P}^m$.
\end{abstract}

\maketitle

\begin{section}{Introduction}
Let $X$ be a projective variety defined over an infinite field $k$ and $\phi$ a finite self morphim of $X$. 
We say $\phi$ is {\it polarized} by a line bundle $L$ on $X$ if $\phi^*L \cong L^{\otimes q}$ for $q>1$.  We say that the polarization is very ample if the line bundle $L$ is very ample i.e. the morphism $X\to \mb{P}(H^0(X,L)) \cong \mb{P}^m_k$ obtained by evaluating the sections of $L$ at points of $X$ is a closed embedding (\cite{H}, pp. 151). In this paper we show that there exists an integer $r\geq1$ such that $\phi^r$ extends to a finite self map of $\mb{P}^m_k$, where $X$ is embedded. We give an example where $r>1$ is required. Fakhruddin proves the following in (\cite{F}, Cor 2.2):

Let $X$ be a projective variety over an infinite field $k$, $\phi:X \to X$ a morphism and $L$ an ample line bundle on $X$ such that $\phi^*L \cong L^{\otimes d}$ for some $d\geq 1$. Then there exists an embedding $i$ of $X$ in $\mb{P}^N_k$ and a morphism $\psi: \mb{P}^N_k \to \mb{P}^N_k$ such that $\psi \circ i = i \circ \phi$.

Our proof and Fakhruddin's proof are closely related. We explain the differences and similarities in the proofs among the two papers in the third remark at the end of this article.

{\it Acknowledgements}: We thank Laura DeMarco and Tom Tucker for their suggestions in the preparation of the paper. 

\end{section}

\begin{section}{Main result}

\bt Let $X$ be a projective variety defined over an infinite field $k$, $L$ a very ample line bundle on $X$ and $\phi: X \to X$ a polarized morphism. Then there exists a positive integer $r$ and a finite morphism $\psi: \mb{P}^m_k \to \mb{P}^m_k$ extending $\phi^r$, where $m +1 = \dim_k H^0(X,L)$. 
\et
\bpf Let $\dim(X)=g$ and let $s_0, \ldots, s_m$ be a basis of $H^0(X,L)$. Let $\mc{I}$ be the sheaf of ideals on $\mb{P}^m$ defining $X$. Then 
\beq \label{SES 1}
0 \to \mc{I} \to \mc{O}_{\mb{P}^m} \to \mc{O}_X \to 0
\eeq
is a short exact sequence of sheaves on $\mb{P}^m$. Tensoring
(\ref{SES 1}) with $\mc{O}_{\mb{P}^m}(n)$ and taking cohomology we get
the long exact sequence 
\beqs
0 \to H^0(\mb{P}^m, \mc{I}(n)) \to H^0(\mb{P}^m,\mc{O}_{\mb{P}^m}(n))
\to
\eeqs
\beqs
\to H^0(X, L^{\otimes n}) \to H^1(\mb{P}^m,\mc{I}(n)) \to \ldots
\eeqs
By Serre's vanishing theorem there
exists $n_0$ depending on $\mc{I}$ such that
$H^1(\mb{P}^m,\mc{I}(n)) =0$ for each $n\geq n_0$.
Let $\{f_i\}$ be the set of homogeneous polynomials defining $X$. Choose
an integer $r$ such that $q^r > \max_i\{\deg(f_i), n_0\}$.
Since $(\phi^r)^*L \cong L^{\otimes q^r}, (\phi^r)^*(s_i)$ can be
lifted to a homogeneous polynomial $h_i$ of degree $q^r$ in the $s_i$'s
defined up to an element of $H^0(\mb{P}^m,\mc{I}(q^r))$. The
polynomials $h_i, 0\leq i \leq m$ define a rational map $\psi: \mb{P}^m \dashrightarrow \mb{P}^m$. We show using induction that if the $h_i$'s
are chosen appropriately, then $\psi$ is a morphism.

Let $W_i$ be the hypersurface defined by $h_i$. We can choose $s_0, \ldots, s_g$ with no common zeros on $X$, then each component (say $Z$) of $\cap_{i=0}^g W_i$
has codimension at most $g+1$ since it is defined by $g+1$ equations. By (\cite{H}, Thm 7.2, pp. 48), it
follows that codim($Z$) $\geq g+1$. Suppose we have $h_0, \ldots, h_j, 0 \leq j \leq m$ such that each component of $\cap_{i=0}^j W_i$ has
codimension $j+1$ and we want to choose $h_{j+1}$. Let $\alpha_1$ be
the lifting of $(\phi^r)^*(s_{j+1})$ to
$H^0(\mb{P}^m,\mc{O}_{\mb{P}^m}(q^r))$. If $V(\alpha_1)$ does not
contain any of the components of $\cap_{i=0}^j W_i$, then set
$h_{j+1} = \alpha_1$.
Otherwise we invoke the  Prime Avoidance Lemma which states:
\bl
Let $A$ be a ring and let $\mf{p}_1, \ldots, \mf{p}_m, \mf{q}$ be ideals of $A$.
Suppose that all but possibly two of the $\mf{p}_i$'s
are prime ideals. If $\mf{q} \nsubseteq \mf{p}_i$ for each $i$, then
$\mf{q}$ is not contained in the set theoretical union $\cup \mf{p}_i$.
\el
\bpf \cite{M}, pp. 2. 
\epf

Taking $A= k[s_0,\ldots, s_n], \mf{q} = \mc{I}(q^r)$, and $\mf{p}_i$'s
the ideals corresponding to the distinct components of $\cap_{i=0}^j W_i$ we can choose
$\alpha_2 \in H^0(\mb{P}^m,\mc{I}(q^r))$ such that $V(\alpha_2)$ does not contain any of
the components of $\cap_{i=0}^j W_i$. Consider the family of hypersurfaces
$V(a\alpha_1 +b\alpha_2)$ with $[a:b] \in \mb{P}^1_k$. If $a=0$, then the corresponding
hypersurface does not contain any components of $\cap_{i=0}^j W_i$. Otherwise, since
$k$ is infinite there exists $c\in k$ such that $V(\alpha_1 + c \alpha_2)$
does not contain any component of $\cap_{i=0}^j W_i$. Let
$h_{j+1} = \alpha_1+c \alpha_2$. This concludes the induction and the
proof of the theorem. 
\epf

We give an example of a self map of a rational quintic in $\mb{P}^3$ that does not extend to $\mb{P}^3$. This illustrates that the condition $r>1$ in Theorem 1 is at times necessary and answers a question of L. DeMarco \cite{D}.

\bp Let $u,v$ be the coordinates of $C\cong \mb{P}^1$ embedded in $\mb{P}^3$
with coordinates $(x_0 = u^5, x_1 = u^4 v, x_2 = uv^4, x_3 = v^5)$.
Then a self map $\phi$ of $C$ of degree $2$ defined by two homogeneous
polynomials $P(u,v)$ and $Q(u,v)$ does not extend to $\mb{P}^3$ if
$P(u,v) = au^2 + buv + cv^2$ with $abc \neq 0$.
\ep
\bpf Considering the restriction map $H^0(\mb{P}^3, \mc{O}_{\mb{P}^3}(2)) \to H^0(\mb{P}^1, \mc{O}_{\mb{P}^1}(10))$. The image of $x^2_0, x^2_1, x^2_2, x^2_3, x_0x_1, x_0x_2, x_0x_3, x_1x_2, x_1x_3, x_2x_3$ under this map is $u^{10}, u^8v^2, u^2v^8, v^{10}, u^9v, u^6v^4, u^5v^5, u^4v^6, uv^9$. Thus $u^7v^3$ and $u^3v^7$ are linearly independent. (Note that it is easy to find two quadratic equations for $C$). One has the possible commutative diagram:
$$ \xymatrix{ \mb{P}^1 \ar[r]^{\phi} \ar[d]^i & \mb{P}^1 \ar[d]^i \\
\mb{P}^3 \ar@{-->}[r]^{\psi} & \mb{P}^3 } $$

The composition ($i \circ \phi$) %with the embedding of $C$ in $\mb{P}^3$ 
is given by four homogeneous polynomials of degree 10, namely $(P(u,v)^5, P(u,v)^4Q(u,v), P(u,v)Q(u,v)^4, Q(u,v)^5)$. If $\phi$ extended to a self map $\psi$ of $\mb{P}^3$, some degree two homogeneous polynomial $F_i$ in the $x_i$'s will restrict to $(P(u,v)^5, P(u,v)^4Q(u,v), P(u,v)Q(u,v)^4, Q(u,v)^5)$ on $C$, by substituting the expressions of the $x_i$ in $(u,v)$. Since $abc \neq 0$ the coefficients of $u^7v^3$ and $u^3v^7$ in $P(u,v)^5$ are non-zero. So $P(u,v)^5$ is not in the image of
$H^0(\mb{P}^3,\mc{O}_{\mb{P}^3}(2)) \to H^0(\mb{P}^1, \mc{O}_{\mb{P}^1}(10))$.
\epf

\end{section}

\begin{section}{Remarks}
\ben
\item If $k$ is finite, $\phi^r$ extends to $\psi$ if we allow $\psi$ to be defined over a finite extension of $k$. Indeed, applying the theorem to $\bar{k}$(the algebraic closure of $k$), $\psi$ is defined by $m+1$ polynomials in $m+1$ variables with coefficients in $\bar{k}$. Hence $\psi$ is defined over the finite extension of $k$ containing the finite set of coefficients of these polynomials. 

\item We say $P\in X$ is {\it preperiodic} for $\phi$ if $\phi^m(P) = \phi^n(P)$ for $m>n\geq 1$. Denote the set of preperiodic points of the dynamical system $(X,L,\phi)$ by Prep($\phi$). It can be easily verified that Prep($\phi$) = Prep($\phi^r$). Thus from an algebraic dynamics perspective, we do not lose any information by replacing $\phi$ by $\phi^r$. The same holds true for points of canonical height \cite{CS} zero as well.

\item One of the technical conditions required to extend $\phi$ from a self map of $X$ to a self map of $\mb{P}^m_k$ is that $\phi^*L \cong L^{s}$ where $s$ is larger than the degrees of equations defining $X$. We choose to replace $\phi$ by $\phi^r$ and fix $L$. The integer $q$ being at least $2$ gives the result immediately. On the other hand in (\cite{F}, Prop 2.1) Fakhruddin chooses to replace $L$ by $L^{\otimes n}$ and keeps $\phi$ fixed. To finish the proof he uses a result of Castelnuevo-Mumford (\cite{Mu}, Theorem 1 and 3), stating that if $n$ is large enough, $X$ will be defined by equations of degree at most two in $\mb{P}(H^0(X,L^{\otimes n}))$.

\een

\end{section}

\end{document}